\pdfoutput=1
\documentclass[twoside]{article}

\usepackage{lipsum} 

\usepackage[sc]{mathpazo} 
\usepackage[T1]{fontenc} 
\usepackage{microtype} 

\usepackage[hmarginratio=1:1,top=32mm,columnsep=20pt]{geometry} 
\usepackage{multicol} 
\usepackage[hang, small,labelfont=bf,up,textfont=it,up]{caption} 
\usepackage{booktabs} 

\usepackage[bookmarks=false]{hyperref} 

\usepackage{lettrine} 
\usepackage{paralist} 

\usepackage{abstract} 

\usepackage{titlesec} 
\titleformat{\section}[block]{\large\scshape}{\thesection.}{1em}{} 
\titleformat{\subsection}[block]{\large}{\thesubsection.}{1em}{} 

\usepackage{fancyhdr} 
\usepackage{amsthm}
\usepackage{amsmath}
\usepackage{amstext}
\usepackage{graphicx}
\usepackage[position=t,singlelinecheck=off]{subfig}
\usepackage{multirow}
\usepackage{tikz}
\usepackage{floatrow}
\usepackage{tipa}
\usepackage{caption}

\newtheorem{dummy}{***}
\newtheorem{thm}[dummy]{Theorem}
\newtheorem{defn}[dummy]{Definition}
\newtheorem{rmk}[dummy]{Remark}

\newtheorem{notn}[dummy]{Notation}
\newtheorem{lem}[dummy]{Lemma}

\usepackage{parskip}
\usepackage{lipsum}
\makeatletter
\def\thm@space@setup{%
  \thm@preskip=12pt plus 0pt minus 8pt
}
\makeatother

\pgfmathsetmacro{\minsize}{0.2cm}

\title{\vspace{-15mm}\fontsize{16pt}{10pt}\selectfont\textbf{A Classification of Non-Compact Coxeter Polytopes with $n+3$ Facets and One Non-Simple Vertex}} 

\author{
\large
\textsc{Mike Roberts}
\\[2mm] 
\normalsize Liverpool University \\ 
\normalsize \href{mailto:michael.roberts@liverpool.ac.uk}{michael.roberts@liverpool.ac.uk} 
\vspace{-5mm}
}
\date{}

\begin{document}

\maketitle 

\thispagestyle{fancy} 


\begin{abstract}

\noindent 
In this paper we state a full classification for Coxeter polytopes in $\mathbb{H}^{n}$ with $n+3$ facets which are non-compact and have precisely one non-simple vertex.

\end{abstract}


\section{Introduction}

A polytope $P$ in hyperbolic $n$-space $\mathbb{H}^{n}$ is called a Coxeter polytope if bounded by hyperplanes which intersect at angle $\frac{\pi}{m_{ij}}$, $m_{ij}\in\mathbb{N}, m_{ij}\ge 2$ for hyperplanes $H_{i}$ and $H_{j}$. They have the interesting property that they tessellate the space $\mathbb{H}^{n}$ through reflections in the bounding hyperplanes.

We can define $P$ as $P=\cap_{i\in J} H_{i}^{-}$, $J$ an arbitrary index set and $H_{i}^{-}$ the half-space containing  $P$. The intersection of a Coxeter polytope with a bounding hyperplane is called a facet of the polytope.

If the Coxeter polytope $P$ has a vertex at infinity, i.e. the hyperplanes intersect at infinity to form a vertex, then $P$ is said to be non-compact, otherwise we say $P$ is compact. Also, if all vertices of $P\subset\mathbb{H}^{n}$ are formed by the intersection of precisely $n$ hyperplanes then $P$ is said to be simple, otherwise we say $P$ is non-simple.

All Coxeter polytopes with $n+1$ facets (simplices) have been fully classified \cite{chein1969, coxeter1950world, lannr1950complexes}. Coxeter polytopes with $n+2$ facets have been fully classified by Kaplinskaya \cite{kaplinskaya1974discrete} for simplicial prisms, Esselmann \cite{Esselmann:1996aa} for compact Coxeter polytopes and Tumarkin \cite{Tumarkin:2004ab} for non-compact Coxeter polytopes.

In the case of Coxeter polytopes with $n+3$ facets Esselmann \cite{esselmann1994kompakte} began and Tumarkin \cite{Tumarkin:2007aa} completed the classification in the compact case. The non-compact case remains without classification currently, however there are some results.

In Tumarkin \cite{Tumarkin:2004aa} there is the following theorem:

\begin{thm}[Tumarkin \cite{Tumarkin:2004aa} Theorem 1] There are no Coxeter polytopes of finite volume with $n+3$ facets in $\mathbb{H}^{n}$ of dimension $n\ge 17$. There is just one such polytope in $\mathbb{H}^{16}$ and it has the following Coxeter diagram:
\begin{figure}[h]
\includegraphics[scale=0.3]{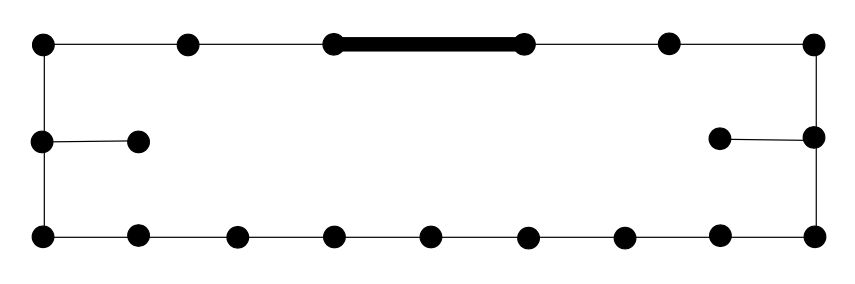}
\end{figure}
\end{thm}

Andreev \cite{andreev1970convex} proved that in dimension $3$ there are infinitely many Coxeter polytopes with $n+3$ facets.

Non-simple pyramids have been classified by Tumarkin \cite{Tumarkin:2004aa} and will not be considered in this text. Therefore a classification is desired for dimensions $4$ to $16$ for non-compact non-pyramidal Coxeter polytopes with $n+3$ facets.

This paper focusses on those polytopes with precisely one non-simple vertex, i.e. the Coxeter polytopes in $\mathbb{H}^{n}$ which have one vertex formed by the intersection of $n+1$ hyperplanes. A full classification of these is displayed in Appendix~\ref{sec:newpoly}.

The majority of this text is copied directly from my Masters dissertation in which a near classification of the non-compact non-pyramidal Coxeter polytopes with $n+3$ facets and one non-simple vertex was obtained.


\section{Gram matrices}\label{sec:grammatrices}

Gram matrices are an elegant way to display and encode the structural information about a Coxeter polytope in matrix form. They also provide the opportunity to question the structure of the Coxeter polytope from a Linear Algebra approach by, for example, considering the eigenvalues and determinant of the Gram matrix.

\begin{defn}[Vinberg \cite{Vinberg:1985aa} \S 1.3]\label{defn:grammatrix} The Gram matrix $G=(g_{ij})$ of a Coxeter polytope $P=\bigcap_{i\in J}H_{i}^{-}\subset\mathbb{H}^{n}$ bounded by $|J|$ hyperplanes is a $|J|\times |J|$ matrix where:
\begin{enumerate}
\item{$g_{ij}=1$ if $i=j$.}
\item{$g_{ij}=-cos(\frac{\pi}{m_{ij}})$ if $H_{i},H_{j}$ intersect.}
\item{$g_{ij}=-1$ if $H_{i},H_{j}$ are parallel.}
\item{$g_{ij}=-cosh(\rho(H_{i},H_{j}))$ if $H_{i}$ and $H_{j}$ are ultraparallel, where $\rho(H_{i},H_{j})$ is the distance between $H_{i}$ and $H_{j}$.}
\end{enumerate}
\end{defn}

Constructing a Gram matrix in this way encodes the angles or distances between the hyperplanes that bound $P$ into a matrix. The matrix is square and from its construction it can be seen to be symmetric, with the number of columns and rows equal to the number of hyperplanes bounding $P$.

\begin{rmk}\label{rmk:detzero}
Consider the columns of the Gram matrix as vectors. $\mathbb{H}^{n}$ can be represented as $\mathbb{E}^{n,1}$ ($n+1$ dimensional Euclidean space with Lorentzian metric). If there are more than $n+1$ columns (vectors) there is a linear dependence between some of the rows or columns and so the Gram matrix has determinant zero. 

This text will consider Coxeter polytopes in $\mathbb{H}^{n}$ with $n+3$ facets and so the Gram matrix of the Coxeter polytope will have determinant zero. A minor formed from $n+2$ rows and columns of this Gram matrix will similarly have determinant zero.

\end{rmk}


\section{Coxeter diagrams \label{sec:coxdiagrams}}

Coxeter diagrams encode into a two-dimensional diagram the structure of the bounding hyperplanes of the Coxeter polytope. A Coxeter polytope has a unique Coxeter diagram.

\begin{defn}[Vinberg \cite{Vinberg:1985aa} $\S 5.1$]\label{def:vincox}For a Coxeter polytope $P$, construct a graph where the vertices $a_{i}$ represent the hyperplanes $H_{i}$ and each node is joined to every other by either a solid, bold or dashed line. The edge between $a_{i}, a_{j}$ is one of the following:
\begin{itemize}
\item{$(m_{ij}-2)-$fold line which is unlabelled (or single line labelled $m_{ij}$) where the dihedral angle between $H_{i}$ and $H_{j}$ is $\frac{\pi}{m_{ij}}$ (note that when $m_{ij}=2$ this is represented with no line),}
\item{a thick line (or single line labelled $\infty$) when the hyperplanes $H_{i}$ and $H_{j}$ are parallel; or}
\item{a dashed line labelled $cosh(\rho(H_{i},H_{j}))$ when the hyperplanes $H_{i}$ and $H_{j}$ are ultraparallel, where $\rho(H_{i},H_{j})$ is the distance between $H_{i}$ and $H_{j}$.}
\end{itemize}
\end{defn}

By reference to Definition~\ref{defn:grammatrix} we see that the Coxeter diagram for a given polytope has the same information contained in it as in the Gram matrix for that polytope, i.e. each Coxeter diagram corresponds to a unique Gram matrix. Therefore the Coxeter diagram defines the combinatorics of the polytope by Vinberg \cite{Vinberg:1985aa} \S 3.


\section{Gale diagrams \label{sec:galediagrams}}

A Gale diagram is obtained by Gale Transform (linear algebra operation) on the vertices of the polytope. Gale transforms will not be discussed here but for more information see, for example, Gr\"{u}nbaum \cite{Grunbaum:1967aa}.

A Gale diagram is of the form $\mathbb{S}^{p-n-2}$ where $p$ is the number of bounding hyperplanes and $n$ is the dimension of the space. There is a node (vertex) labelled $1$ on the boundary of the $(p-n-2)$-sphere for each of the $k$ hyperplanes.

This text focusses on Coxeter polytopes with $n+3$ bounding hyperplanes and hence the Gale diagram is $\mathbb{S}^{1}$, i.e. a circle, and the nodes for the hyperplanes are on the circumference of it.

We use a reformulation of the definition of the Gale diagram taken from Tumarkin \cite{Tumarkin:2004aa} p1 which will give a standard (contracted) form for Gale diagrams.

\begin{defn}[Gr\"{u}nbaum \cite{Grunbaum:1967aa}]\label{defn:galediagrams}
Every combinatorial type of an $n$-dimensional polytope with $n + 3$ facets can be represented by a standard two-dimensional Gale diagram. This consists of vertices of regular $2k-$gon in $\mathbb{E}^{2}$ centered at the origin which are labelled according to the following rules:
\begin{enumerate}
\item{Each label is a non-negative integer, the sum of labels equals $n + 3$.}
\item{Labels of neighbouring vertices cannot be equal to zero simultaneously.}
\item{Labels of opposite vertices can not be equal to zero simultaneously.}
\item{The points that lie in any open halfspace bounded by a hyperplane through the origin have labels whose sum is at least two.}
\end{enumerate}
\end{defn}

\begin{rmk}[Tumarkin \cite{Tumarkin:2004aa} p1]\label{rmk:galefaces}
The combinatorial type of a convex polytope can be read off from the Gale diagram in the following way. Each vertex $a_{i}, i = 1, . . . , 2k$, with label $\mu_{i}$ corresponds to $\mu_{i}$ facets $f_{i,1},\ldots,f_{i,\mu_{i}}$ of $P$. For any subset $I$ of the set of facets of $P$ the intersection of facets $\{ f_{j,\gamma} | (j,\gamma )\in I \}$ is a face of $P$ if and only if the origin is contained in the set $conv \{ a_{j} | ( j , \gamma ) \not\in I \}$.
\end{rmk}


\section{Methodical approach for constructing Coxeter polytopes}

\subsection{Introduction}

In this section we will detail the rigorous approach taken to classify the non-compact non-pyramidal Coxeter polytopes with $n+3$ facets and one non-simple vertex. We start by introducing some notation and then two important lemmas.

\begin{notn}[Tumarkin \cite{Tumarkin:2004aa} \S 2 p237] If $G$ is the Gale diagram of a polytope $P$, denote by $S_{m,l}$ the subdiagram of the Coxeter diagram $S(P)$ corresponding to the $l-m+1$ (mod $2k$) consecutive vertices $a_{m},\ldots, a_{l}$ of $G$. For $l=m$, denote this subdiagram by $S_{m}$. The weight of the vertex is denoted by $\mu(a_{i})$.
\end{notn}

\begin{lem}[Tumarkin \cite{Tumarkin:2004aa} Lemma 1]\label{lem:tum1} Let $G$ be the Gale diagram of a polytope $P$. Suppose that the weights of $a_{i},a_{k+i}$ are non-zero. Then:
\begin{enumerate}
\item{the vertices $a_{i}$ and $a_{k+i}$ have weight $1$ and the Coxeter diagrams $S_{i+1,k+i-1}$ and $S_{k+i+1,i-1}$ are connected and parabolic.}
\item{if $a_{i+1}$ and $a_{k+i+1}$ have non-zero weights, the Coxeter diagram $S_{i+1,k+i}$ is quasi-Lann\'{e}r.}
\item{if $a_{i+1}$ has weight zero, the Coxeter diagram $S_{i+2,k+i}$ is quasi-Lann\'{e}r.}
\end{enumerate}
\end{lem}

\begin{lem}[Tumarkin \cite{Tumarkin:2004aa} Lemma 2] \label{lem:tum2} Let $G$ be the Gale diagram of a polytope $P$, and suppose that the weights of the vertices $a_{i}$ and $a_{k+i-1}$ are zero. Then the diagram $S_{i+1,k+i-2}$ is Lann\'{e}r.
\end{lem}

\begin{rmk}We see that a non-simple Coxeter polytope is represented by a Gale diagram which has pairs of nodes opposite one another around the circumference and whose labels must be $1$ (Lemma~\ref{lem:tum1}.1). In this text we focus on those Gale diagrams which are a circle with precisely one pair of nodes opposite one another.
\end{rmk}

Notice that there are only a finite number of Gale diagrams in each dimension which have one pair of nodes opposite each other, this can be seen as the task of splitting integer $n+1$ over the nodes which are not opposite one another.

As the final point before detailing the methodical approach taken, we quote the following important theorem:

\begin{thm}[Vinberg \cite{Vinberg:1985aa} Theorem 2.1]\label{thm:vin21} Let $G=(g_{ij})$ be an indecomposable symmetric matrix of signature $(n,1)$ with $1's$ along the diagonal and non-positive entries off it. Then there is a convex polytope $P$ in $\mathbb{H}^{n}$ whose Gram matrix is $G$. The polytope $P$ is uniquely determined up to isometry in $\mathbb{H}^{n}$.
\end{thm}

\subsection{The method explained}

Suppose we look to find all non-compact non-pyramidal Coxeter polytopes in $\mathbb{H}^{n}$ which have $n+3$ facets and one non-simple vertex, then applying the following method will give all such Coxeter polytopes.

\begin{enumerate}
\item{Determine all possible Gale diagrams up to congruence.}
\item{Apply Tumarkin's Lemma's~\ref{lem:tum1} and \ref{lem:tum2}. These will restrict Coxeter subdiagrams which can be formed.}
\item{Apply the restrictions which can be read off from the Gale diagram such as those subdiagrams which are necessarily elliptic, etc.}
\item{Form a Gram matrix $G$ for the candidate Coxeter diagrams and check that the determinant of $G$ and of any $(n+2)\times(n+2)$ minor is zero (Remark~\ref{rmk:detzero})}
\item{Apply Vinberg's Theorem~\ref{thm:vin21} and check the signature of the Gram matrix $G$. If the signature is $(n,1,2)$ then the Gram matrix $G$ corresponds to a Coxeter polytope.}
\end{enumerate}


\section{Results and Further Work}

The full classification for the non-compact non-pyramidal Coxeter polytopes with $n+3$ facets and precisely one non-simple vertex is shown in Appendix~\ref{sec:newpoly}. It is interesting to note that there are no examples in dimension 11 and above.

To obtain a full classification of the non-compact Coxeter polytopes with $n+3$ facets it remains to classify such simple Coxeter polytopes and non-simple non-pyramidal Coxeter polytopes with more than one non-simple vertex. 

Once this is completed then along with along with \cite{Tumarkin:2007aa} and \cite{Tumarkin:2004aa} this will complete the classification of Coxeter polytopes with $n+3$ facets.

\vspace{0.2in}
{\bf Acknowledgements}

The majority of this work was completed at Durham University as a part of my Masters thesis. I am extremely grateful to my supervisor Dr Pavel Tumarkin for all his help and advice whilst I was performing the work and for his helpful comments on this paper. I'm also thankful to Rafael Guglielmetti for his efforts in verifying the Coxeter polytopes found and providing me with corrections.


\bibliography{biblio}

\begin{thebibliography}{10}

\bibitem{andreev1970convex}
E.M. Andreev.
\newblock On convex polyhedra of finite volume in {L}obachevskii spaces.
\newblock {\em USSR Sbornik}, 12:255--259, 1970.

\bibitem{chein1969}
M.~Chein.
\newblock Recherche des graphes des matrices de {C}oxeter hyperboliques d'ordre
  $\le 10$.
\newblock {\em Revue fran{\c c}aise d'informatique et de recherche
  op{\'e}rationnelle}, 3:3--16, 1969.

\bibitem{coxeter1950world}
H.S.M. Coxeter and G.J. Whitrow.
\newblock World-structure and non-euclidean honeycombs.
\newblock In {\em Proceedings of the Royal Society of London A: Mathematical,
  Physical and Engineering Sciences}, volume 201, pages 417--437. The Royal
  Society, 1950.

\bibitem{esselmann1994kompakte}
F.~Esselmann.
\newblock {\"U}ber kompakte, hyperbolische coxeter-polytope mit wenigen
  facetten.
\newblock {\em Preprint/Universit{\"a}t Bielefeld, Sonderforschungsbereich 343,
  Diskrete Strukturen in der Mathematik; 94-087}, 1994.

\bibitem{Esselmann:1996aa}
F.~Esselmann.
\newblock The classification of compact hyperbolic {C}oxeter $d-$polytopes with
  $d+2$ facets.
\newblock {\em Commentarii Mathematici Helvetici}, 71(1):229--242, 1996.

\bibitem{Grunbaum:1967aa}
B.~Grunba\"{u}m, V.~Klee, M.A. Perles, and G.C. Shephard.
\newblock {\em Convex polytopes}.
\newblock Springer, 1967.

\bibitem{kaplinskaya1974discrete}
I.M. Kaplinskaya.
\newblock Discrete groups generated by reflections in the faces of symplicial
  prisms in {L}obachevskian spaces.
\newblock {\em Mathematical Notes of the Academy of Sciences of the USSR},
  15(1):88--91, 1974.

\bibitem{lannr1950complexes}
F.~Lann{\'e}r.
\newblock On complexes with transitive groups of automorphisms.
\newblock {\em Comm. Sem. Math. Univ. Lund}, 11:1--71, 1950.

\bibitem{Tumarkin:2004ab}
P.~Tumarkin.
\newblock Hyperbolic {C}oxeter n-polytopes with $n+2$ facets.
\newblock {\em Mathematical Notes}, 75(5-6):848--854, 2004.

\bibitem{Tumarkin:2004aa}
P.~Tumarkin.
\newblock Hyperbolic {C}oxeter n-polytopes with $n+3$ facets.
\newblock {\em Transactions of the Moscow Mathematical Society}, 65:235--250,
  2004.

\bibitem{Tumarkin:2007aa}
P.~Tumarkin.
\newblock Compact hyperbolic {C}oxeter $n-$polytopes with $n+3$ facets.
\newblock {\em Journal of Combinatorics}, 14(4):R69, 2007.

\bibitem{Vinberg:1985aa}
E.B. Vinberg.
\newblock Hyperbolic reflection groups.
\newblock {\em Russian Mathematical Surveys}, 40(1):31--75, 1985.

\end{thebibliography}
\bibliographystyle{plain}


\clearpage
\appendix
\section{New Coxeter Polytopes \label{sec:newpoly}}

This section details the new Coxeter polytopes found which have $n+3$ facets, are non-compact non-pyramidal and have precisely one non-simple vertex. They are split over the following subsections by the dimension in which they exist. 

The nodes are labelled by the location of that node in its Gale diagram. For example, a node labelled $j$ will appear at location $a_{j}$ in the Gale diagram with label $1$, a node with label $i,j$ (with $j\in\{ 1,\ldots s\}$) will be one hyperplane from node $a_{i}$ in the Gale diagram (which has label $s$). This is purely to provide an easy trail back to the Gale diagram for the Coxeter polytope and these node labels may be ignored.

\clearpage
\subsection{Dimension 4}

\begin{figure}[h]

\sidesubfloat[]{%

}
\end{figure}

\clearpage

\subsection{Dimension 5}

\begin{figure}[!htb]
\sidesubfloat[]{%
%
}

\end{figure}

\clearpage
\subsection{Dimension 6}

\begin{figure}[!htb]
\sidesubfloat[]{%
%
}

\end{figure}
\clearpage
\subsection{Dimension 7}

\begin{figure}[!htb]
\sidesubfloat[]{%
%
}

\end{figure}
\clearpage

\subsection{Dimension 8}

\begin{figure}[!htb]
\sidesubfloat[]{%
%
}

\end{figure}
\clearpage
\subsection{Dimension 9}

\begin{figure}[!htb]
\sidesubfloat[]{%
%
}

\end{figure}
\clearpage
\subsection{Dimension 10}

\begin{figure}[!htb]
\sidesubfloat[]{%
%
}
\end{figure}

\end{document}